\documentclass[11pt]{article}

\usepackage{amsmath,amsthm,amsfonts}

\usepackage{amssymb}

\usepackage{epsfig}

\usepackage{rotating}

\usepackage{subfigure}

\newcommand{\p}[2]{\frac{\partial#1}{\partial#2}}

\def\hR{{\hat{R}}}

\begin{document}

\title{Phase transition approach to detecting singularities of PDEs}
\author{Panagiotis Stinis \\ 
Department of Mathematics \\
University of Minnesota \\
    Minneapolis, MN 55455} 
\date {}

\maketitle

\begin{abstract}
We present a mesh refinement algorithm for detecting singularities of time-dependent partial differential equations. The algorithm is inspired by renormalization constructions used in statistical mechanics to evaluate the properties of a system near a critical point, i.e. a phase transition. The main idea behind the algorithm is to treat the occurrence of singularities of time-dependent partial differential equations as phase transitions. 

The algorithm assumes the knowledge of an accurate reduced model. In particular, we need only assume that we know the functional form of the reduced model, i.e. the terms appearing in the reduced model, but not necessarily their coefficients. We provide a way of computing the necessary coefficients on the fly as needed. 

We show how the mesh refinement algorithm can be used to calculate the blow-up rate as we approach the singularity. This calculation can be done in three different ways: i) the direct approach where one monitors the blowing-up quantity as it approaches the singularity and uses the data to calculate the blow-up rate ; ii) the ``phase transition" approach (\`a la Wilson) where one treats the singularity as a fixed point of the renormalization flow equation and proceeds to compute the blow-up rate via an analysis in the vicinity of the fixed point and iii) the ``scaling" approach (\`a la Widom-Kadanoff) where one postulates the existence of scaling laws for different quantities close to the singularity, computes the associated exponents and then uses them to estimate the blow-up rate. Our algorithm allows a unified presentation of these three approaches. 

The inviscid Burgers and the supercritical focusing Schr\"odinger equations are used as instructive examples to illustrate the constructions.  

\end{abstract}


\section*{Introduction}

The problem of how to construct mesh refinement methods and how to approach more efficiently possible singularities of partial differential equations has attracted considerable attention (see e.g. \cite{almgren,berger1,berger2,budd,ceniceros,landman}). At the same time, the problem of constructing dimensionally reduced models for large systems of ordinary differential equations (this covers the case of partial differential equations after discretization or series expansion of the solution) has also received considerable attention e.g. see the review papers \cite{givon, CS05}. The construction of an accurate reduced model has advantages beyond the obvious one of predicting the correct behavior for a reduced set of variables.

We present here an algorithm that is based on dimensional reduction and which can be used to perform mesh refinement and investigate possibly singular solutions of partial differential equations (see also \cite{S07}). The algorithm is inspired by constructions used in statistical mechanics to evaluate the properties of a system near a critical point \cite{goldenfeld, binney} (a critical point is a value for the controlling parameter of a system at which the behavior of the system changes abruptly). The idea underlying the computation of the properties at criticality,  is that while the form of the reduced system equations is important, one can extract even more information by looking at how the form of the reduced system is related to the form of the original (full dimensional) system \cite{weinberg, wilson}. We extend this idea to the study of (possibly) singular solutions of partial differential equations by treating time as the controlling parameter and the instant of occurrence of a singularity as a critical value for the parameter, i.e. a critical point.

Our approach has two objectives: i) it provides a way of accurately monitoring the progress of a simulation towards underresolution, thus providing us as a byproduct with the time of occurrence of the possible singularity; ii) it allows the formulation of a mesh refinement scheme which is able to reach the time of interesting dynamics of the equation much more efficiently compared to an algorithm that simply starts with the maximum available resolution. 

The mesh refinement algorithm can be used to calculate the blow-up rate as we approach the singularity. This calculation can be done in three different ways: i) the direct approach where one monitors the blowing-up quantity as it approaches the singularity and uses the data to calculate the blow-up rate ; ii) the ``phase transition" approach (\`a la Wilson) \cite{goldenfeld} where one treats the singularity as a fixed point of the renormalization flow equation and proceeds to compute the blow-up rate via an analysis in the vicinity of the fixed point and iii) the ``scaling" approach (\`a la Widom-Kadanoff) \cite{binney} where one postulates the existence of scaling laws for different quantities close to the singularity, computes the associated exponents and then uses them to estimate the blow-up rate. Our algorithm allows a unified presentation of these three approaches.

The task of investigating numerically the appearance of a singularity is subtle. Clearly, since all calculations are performed with finite resolution and a singularity involves an infinity of active scales we can only come as close to the singularity as our resolution will allow. After that point, either we stop our calculations and conclude that a singularity may be present close to the time instant we stopped or we switch, if available, to a model that drains energy at the correct rate out of the set of resolved variables. We emphasize here that, up to some time, the evolution towards a near-singular solution can be identical to the evolution towards a singular solution. If we do not have enough resolution to go beyond the time instant after which the two evolutions start deviating, we cannot claim with certainty the presence of a singularity. In other words, given adequate resolution we can eliminate the possibility of a singularity. But, it may be very hard, to prove through a finite calculation, that a singularity exists (we come back to these points in Sections \ref{know2} and \ref{numerics2}).

The paper is organized as follows. In Section \ref{algo} we present the ideas behind the construction of the algorithm.  In Section \ref{know2} we present the mesh refinement algorithm. In Section \ref{numerics2} we provide numerical results for the inviscid Burgers equation. In Section \ref{schrodinger} we provide numerical results for the supercritical focusing Schr\"odinger equation. Section \ref{rate} shows how one can use the mesh refinement algorithm to compute the blow-up rate as a critical exponent, i.e. using solely properties of a renormalization (coarse-graining) process in the vicinity of the singularity. Section \ref{conclusions} contains a discussion of the results and some directions for future work.


\section{The main construction}\label{algo}

Suppose that we are interested in the possible development of singularities in the solution $v(x,t)$ of a partial differential equation (PDE)
$$ v_t + H (t,x,v,v_x,...)=0 $$
where $H$ is a, in general nonlinear, operator and $x \in D \subseteq \mathbb{R}^d$ (the constructions extend readily to the case of systems of partial differential equations). After spatial discretization or expansion of the solution in series, the PDE transforms into a system of ordinary differential equations (ODEs). For simplicity we restrict ourselves to the case of periodic boundary conditions, so that a Fourier expansion of the solution leads to system of ODEs for the Fourier coefficients. To simulate the system for the Fourier coefficients we need to truncate at some point the Fourier 
expansion. Let $F \cup G$ denote the set of Fourier modes retained in the series, where we have split the Fourier modes in two sets, $F$ and $G.$ 
We call the modes in $F$ resolved and the modes in $G$ unresolved. One can construct, in principle, an 
exact reduced model for the modes in $F$ e.g. through the Mori-Zwanzig formalism \cite{CHK00} (we do not deal here with the complications of constructing a good reduced model).

The main idea behind the algorithm is that the evolution of moments of the reduced set of modes, for example $l_p$ norms of the modes in $F$, should be the same whether computed from the full or the reduced system. This is a generalization to time-dependent systems of the principle used in the theory of equilibrium phase transitions to compute the critical exponents \cite{goldenfeld,S05}. The idea underlying the computation of the critical exponents is that while the form of the reduced system equations is important, one can extract even more information by looking at how the form of the reduced system is related to the form of the original (full dimensional) system. We extend this idea to the study of (possibly) singular solutions of partial differential equations by treating time as the controlling parameter and the instant of occurrence of a singularity as a critical value for the parameter, i.e. a critical point. We caution the reader that even though our motivation for the present construction came from the theory of equilibrium phase transitions, we do not advocate that a singularity is a phase transition in the conventional sense. It can be thought of as a transition from a strong solution to an appropriately defined weak solution but one does not have to push the analogy further. We want to point here that the problem we are addressing is different from the subject known as dynamic critical phenomena (see Ch. 8 in \cite{goldenfeld}). There, one is interested in the computation of time-dependent quantities as a controlling parameter, other than time, reaches its critical value. In our case, time {\it is} the controlling parameter and we are interested in the behavior of the solution as time reaches a critical value.

The above arguments can be made more precise. The original system of equations for the modes $F \cup G$ is given by 
$$\frac{du(t)}{dt} = R (t,u(t)),$$
where $u = ( \{u_k\}), \; k \in F \cup G$ is the vector of Fourier coefficients of $u$ and $R$ is the Fourier transform of the operator $H.$ The system should be supplemented with an initial condition $u(0)=u_0.$ The vector of Fourier coefficients can be written as $ u = (\hat{u}, 
\tilde{u}),$ where $ \hat{u}$ are the resolved modes (those in $F$) and $\tilde{u}$ the unresolved ones (those in $G$). Similarly, for the right hand sides (RHS) we have $R(t,u) = (\hat{R}(t,u), \tilde{R}(t,u)).$ Note that the RHS of the resolved modes involves both resolved and unresolved modes. In anticipation of the construction of a reduced model we can rewrite the RHS as $R(t,u)=R^{(0)}(t,u) = (\hat{R}^{(0)}(t,u), \tilde{R}^{(0)}(t,u)).$ For each mode $u_k, \; k \in F \cup G,$ we can decompose $R_k^{(0)}(t,u)$ as 
$$R_k^{(0)}(t,u(t)) = \sum_{i=1}^{m} a^{(0)}_i R^{(0)}_{ik}(t,u(t)).$$
Thus, the equation for the the mode $u_k, \; k \in F \cup G$ is written as 
\begin{equation}\label{full}
\frac{d{u_k}(t)}{dt} = {R}_k(t,u)={R}_k^{(0)} (t,u(t))=\sum_{i=1}^{m} a^{(0)}_i {R}^{(0)}_{ik} (t,u(t))
\end{equation}
Note that not all the coefficients $a^{(0)}_i, \; i =1,\ldots,m$ have to be nonzero. As is standard in renormalization theory \cite{binney}, one augments (with zero coefficients) the RHS of the equations in the full system by terms whose form is the same as the terms appearing in the RHS of the equations for the reduced model. Dimensional 
reduction transforms the vector $a^{(0)}=(a^{(0)}_1,\ldots,a^{(0)}_m)$ to $a^{(1)}=(a^{(1)}_1,\ldots,a^{(1)}_m).$ The reduced model for the mode $u'_k, \; k \in F $ is given by 
\begin{equation}\label{reduced}
\frac{d{u}_k'(t)}{dt} = R_k^{(1)} (t,\hat{u}'(t))=\sum_{i=1}^{m} a^{(1)}_i R^{(1)}_{ik} (t,\hat{u}'(t))
\end{equation}
with initial condition ${u}_k'(0)={u}_{0k}.$ We emphasize that the functions $R^{(1)}_{ik}, \; i=1,\ldots,m, \; k \in F,$ have the same form as the functions $R^{(0)}_{ik}, \; i=1,\ldots,m, \; k \in F,$ but they depend {\it only} on the reduced set of modes $F.$ This allows one to determine the relation of the full to the reduced system by focusing on the change of the vector $a^{(0)}$ to $a^{(1)}.$ Also, the vectors $a^{(0)}$ and $a^{(1)}$ do not have to be constant in time. This does not change the analysis that follows.

Define $m$ quantities 
$\hat{E}_i, \; i=1,\ldots,m$ involving only modes in $F.$ For example, these could be $L_p$ norms of the reduced set of modes. To proceed we require that these quantities' rates of change are the same when computed from (\ref{full}) and (\ref{reduced}), i.e. 
\begin{equation}\label{conditions}
\frac{d\hat{E}_i(\hat{u})}{dt} = \frac{d\hat{E}_i(\hat{u}')}{dt}, \; i=1,\ldots,m.
\end{equation}
Note that similar conditions, albeit time-independent, lie at the heart of the renormalization group theory for equilibrium systems (\cite{binney} p. 154). In fact, it is these conditions that allow the definition and calculation of the (renormalization) matrix whose eigenvalues are used to calculate the critical exponents. In the current (time-dependent) setting, the renormalization matrix is defined by 
differentiating $\frac{d\hat{E}_i(\hat{u})}{dt}$ with respect to $a^{(0)}$ and using (\ref{conditions}) to 
obtain
\begin{equation}\label{rngmatrix1}
\p{}{a^{(0)}_j}\biggl(\frac{d\hat{E}_i(\hat{u})}{dt}\biggr)=\sum_{k=1}^{m}\p{}{a^{(1)}_k}\biggl(\frac{d\hat{E}_i(\hat{u}')}{dt}\biggr) \p{a^{(1)}_k}{a^{(0)}_j}, \; i,j=1,\ldots,m.
\end{equation}
We define the renormalization matrix $M_{kj}= \p{a^{(1)}_k}{a^{(0)}_j}, \; k,j=1,\ldots,m,$ as well as the matrices $A_{kj}=\p{}{a^{(0)}_j}\biggl(\frac{d\hat{E}_k(\hat{u})}{dt}\biggr), \; k,j=1,\ldots,m$ and $B_{kj}=\p{}{a^{(1)}_j}\biggl(\frac{d\hat{E}_k(\hat{u}')}{dt}\biggr), \; k,j=1,\ldots,m .$ Equations (\ref{rngmatrix1}) can 
be written in matrix form as
\begin{equation}\label{rngmatrix2}
A=MB
\end{equation} 
The entries of $A$ describe the contributions of the different terms appearing on the RHS of the full system to the rate of change of $E_i$. The same for the entries of matrix $B$ and the reduced model.

The eigenvalues of the matrix $M$ contain information about the behavior of the reduced system relative to the full system. In fact, they measure whether the full and reduced systems deviate or approach. In the renormalization theory of critical phenomena, the eigenvalues of $M$ at the critical point are used to analyze the system properties close to criticality. The analysis is based on the assumption that the eigenvalues of $M$ change slowly near the critical point so that even if one cannot compute exactly {\it on} the critical point, it is possible to get an accurate estimate of them by computations near the critical point. Then, one performs a linear stability analysis near the fixed point and computes the system properties. The situation in the case of singularities of PDEs  is different. In this case, the eigenvalues of $M$ vary {\it most rapidly} near the singularity, due to the full system's rapid deterioration. Thus, we are not able to use linear stability analysis near the singularity. However, we are still able to extract the relevant blow-up rates (see Section \ref{rate}).

\subsection{An instructive example}\label{example}

We use the 1D inviscid Burgers equation as an instructive example for the constructions presented in this section. The equation is given by 
\begin{equation}\label{burgers}
u_t+u u_x = 0.
\end{equation}
Equation (\ref{burgers}) should be supplemented with an initial condition $u(x,0)=u_0(x)$ and boundary conditions. We solve (\ref{burgers}) in the interval $[0,2\pi]$ with periodic boundary conditions. This allows us to expand the solution in Fourier series
$$u^{M}(x,t )=\underset{k \in F \cup G}{\sum} u_k(t) e^{ikx},$$
where $F \cup G=[-\frac{M}{2},\frac{M}{2}-1].$ We have written the set of Fourier modes as the union of two sets 
in anticipation of the construction of the reduced model comprising only of the modes in $F=[-\frac{N}{2},\frac{N}{2}-1],$ where $ N < M.$
The equation of motion for the Fourier mode $u_k$ becomes
\begin{equation}
\label{burgersode}
 \frac{d u_k}{dt}=- \frac{ik}{2} \underset{p, q \in F \cup G}{\underset{p+q=k  }{ \sum}} u_{p} u_{q}.
\end{equation}

\subsubsection{The $t$-model}\label{numericstmodel}

We need to choose a reduced model for the modes in $F.$ We use a reduced model, known as the $t$-model, which follows correctly the behavior of the solution to the inviscid Burgers equation even after the formation of shocks \cite{bernstein,HS06}. The $t$-model was first derived in the context of statistical irreversible mechanics \cite{CHK3} and was later analyzed in \cite{bernstein,HS06}. It is based on the assumption of the absence of time scale separation between the resolved and unresolved modes. We will use the same model for the case with nonzero viscosity and comment on its validity when appropriate. For a mode $u_k'$ in $F$ the model is given by
\begin{multline}\label{burgersode2}
\frac{d}{dt}u_k'=- \frac{ik}{2}   \underset{p \in F ,\, q \in F }{\underset{p+q=k  }{ \sum}}   u_{p}'u_{q}'  \\
-\frac{ik}{2} \underset{p \in F  ,\, q \in G}{\underset{p+q=k  }{ \sum}}   u_{p}'  \biggl[ -t\frac{iq}{2}  \underset{r \in F ,\,  s \in F }{\underset{r+s=q }{ \sum}}  u_{r}'u_{s}' \biggr] \\
-\frac{ik}{2} \underset{p \in G ,\, q \in F }{\underset{p+q=k  }{ \sum}}  \biggl[ -t\frac{ip}{2} \underset{r \in F  ,\, s \in F }{\underset{r+s=p }{ \sum}} u_{r}'
u_{s}' \biggr]  u_{q}' . 
\end{multline}
The first term on the RHS of (\ref{burgersode2}) is of the same form as the first term in (\ref{burgersode}), except that the term in (\ref{burgersode2}) is defined only for the modes in $F.$ The viscous term is the same. The third and fourth terms in (\ref{burgersode2}) are not present in (\ref{burgersode}). They are cubic in the Fourier modes and they are effecting the drain of energy out of the modes in $F.$ We should note here that the cubic terms in the $t$-model do not depend on the viscosity. To conform with the notation in Section \ref{algo} we rewrite (\ref{burgersode2}) as
\begin{multline*}
\frac{d}{dt}u_k'= a^{(1)}_1 \biggl[- \frac{ik}{2}   \underset{p \in F ,\, q \in F }{\underset{p+q=k  }{ \sum}}   u_{p}'u_{q}'   \biggr] +\\
a^{(1)}_2 \Biggl[ -\frac{ik}{2} \underset{p \in F  ,\, q \in G}{\underset{p+q=k  }{ \sum}}   u_{p}'  \biggl[ -t\frac{iq}{2}  \underset{r \in F ,\,  s \in F }{\underset{r+s=q }{ \sum}}  u_{r}'u_{s}' \biggr] \\
-\frac{ik}{2} \underset{p \in G ,\, q \in F }{\underset{p+q=k  }{ \sum}}  \biggl[ -t\frac{ip}{2} \underset{r \in F  ,\, s \in F }{\underset{r+s=p }{ \sum}} u_{r}'
u_{s}' \biggr]  u_{q}' \Biggr],  
\end{multline*}
where $a^{(1)}_1=1$ and $a^{(1)}_2=1.$ We rewrite Equation (\ref{burgersode})  as
\begin{multline*}
 \frac{d u_k}{dt}= a^{(0)}_1 \biggl[ - \frac{ik}{2} \underset{p, q \in F \cup G}{\underset{p+q=k  }{ \sum}} u_{p} u_{q}  \biggr] + \\
a^{(0)}_2 \Biggl[ -\frac{ik}{2} \underset{p \in F\cup G  ,\, q \in I}{\underset{p+q=k  }{ \sum}}   u_{p}  \biggl[ -t\frac{iq}{2}  \underset{r \in F\cup G ,\,  s \in F \cup G }{\underset{r+s=q }{ \sum}}  u_{r}u_{s} \biggr] \\
-\frac{ik}{2} \underset{p \in I ,\, q \in F\cup G }{\underset{p+q=k  }{ \sum}}  \biggl[ -t\frac{ip}{2} \underset{r \in F\cup G  ,\, s \in F\cup G }{\underset{r+s=p }{ \sum}} u_{r} u_{s} \biggr]  u_{q} \Biggr]  ,
\end{multline*}
where $a^{(0)}_1=1$ and $a^{(0)}_2=0.$ The reader should note that we have introduced a new set of modes $I.$ This is the set of unresolved modes for the {\it full} system. The reason for introducing the set $I$ is that, as is the case in renormalization formulations, the terms appearing in the RHS of the equations at the different levels of resolution should be of the same functional form. The difference between the different levels of resolution should be only in the range of modes used. Since the $t$-model involves a quadratic convolution sum with one index in the resolved range and the other in the unresolved range, we should use the same functional form when constructing the corresponding term for the full system. Thus, this term should involve a convolution sum with one index in the range $F \cup G$ and the other in $I.$ 

Further, define
$$\hR^{(0)}_{1k}(t,\hat{u}(t))=- \frac{ik}{2} \underset{p, q \in F \cup G}{\underset{p+q=k  }{ \sum}} u_{p} u_{q} $$
and
\begin{multline*}
\hR^{(0)}_{2k}(t,\hat{u}(t))= -\frac{ik}{2} \underset{p \in F\cup G  ,\, q \in I}{\underset{p+q=k  }{ \sum}}   u_{p}  \biggl[ -t\frac{iq}{2}  \underset{r \in F\cup G ,\,  s \in F \cup G }{\underset{r+s=q }{ \sum}}  u_{r}u_{s} \biggr] \\
-\frac{ik}{2} \underset{p \in I ,\, q \in F\cup G }{\underset{p+q=k  }{ \sum}}  \biggl[ -t\frac{ip}{2} \underset{r \in F\cup G  ,\, s \in F\cup G }{\underset{r+s=p }{ \sum}} u_{r} u_{s} \biggr]  u_{q}
\end{multline*}
Also, define
$$\hR^{(1)}_{1k}(t,\hat{u}'(t))=- \frac{ik}{2} \underset{p, q \in F }
{\underset{p+q=k  }{ \sum}} u_{p}' u_{q}'  $$
and
\begin{multline*}
\hR^{(1)}_{2k}(t,\hat{u}'(t))=-\frac{ik}{2} \underset{p \in F  ,\, q \in G}{\underset{p+q=k  }{ \sum}}   u_{p}'  \biggl[ -t\frac{iq}{2}  \underset{r \in F ,\,  s \in F }{\underset{r+s=q }{ \sum}}  u_{r}'u_{s}' \biggr] \\
-\frac{ik}{2} \underset{p \in G ,\, q \in F }{\underset{p+q=k  }{ \sum}}  \biggl[ -t\frac{ip}{2} \underset{r \in F  ,\, s \in F }{\underset{r+s=p }{ \sum}} u_{r}'
u_{s}' \biggr]  u_{q}'
\end{multline*}
Thus, the equations of motion for the resolved modes in the full system and the reduced model can be written as
\begin{equation}\label{alburgersode}
 \frac{d u_k}{dt}=\sum_{i=1}^{2} a^{(0)}_i \hR^{(0)}_{ik}(t,{u}(t))
\end{equation}
and
\begin{equation}\label{alburgersode2}
\frac{d u_k'}{dt}=\sum_{i=1}^{2} a^{(1)}_i \hR^{(1)}_{ik}(t,\hat{u}'(t))
\end{equation}

To proceed, we need to define the quantities $\hat{E}_i, \; i=1,\ldots,m.$ In our case, $m=2$ and we need to define $\hat{E}_1$ and $\hat{E}_2.$ The choice of the $\hat{E}_i$ is not unique. We chose for our experiments $\hat{E}_1=\underset{k \in F}{\sum} |u_k|^2$ and $\hat{E}_2=\underset{k \in F}{\sum} |u_k|^4.$ The rates of change of the $\hat{E}_i$ are given for the full system by
$$\frac{d\hat{E}_1}{dt}=\underset{k \in F}{\sum}a^{(0)}_1 2 Re(\hR^{(0)}_{1k}(t,\hat{u}(t)) u_k^{*})+ a^{(0)}_2 2 Re(\hR^{(0)}_{2k}(t,\hat{u}(t)) u_k^{*})$$
and
$$\frac{d\hat{E}_2}{dt}=\underset{k \in F}{\sum}a^{(0)}_1 2 Re(2 \hR^{(0)}_{1k}(t,\hat{u}(t)) |u_k|^2 u_k^{*})+ a^{(0)}_2 2 Re(2 \hR^{(0)}_{2k}(t,\hat{u}(t)) |u_k|^2 u_k^{*})$$
where 
$u_k^{*}$ is the complex conjugate of $u_k.$ Similarly, for the reduced system we have
$$\frac{d\hat{E}_1}{dt}=\underset{k \in F}{\sum}a^{(1)}_1 2 Re(\hR^{(1)}_{1k}(t,\hat{u}'(t)) u_k'^{*})+ a^{(1)}_2 2 Re(\hR^{(1)}_{2k}(t,\hat{u}'(t)) u_k'^{*})$$
and
$$\frac{d\hat{E}_2}{dt}=\underset{k \in F}{\sum}a^{(1)}_1 2 Re(2 \hR^{(1)}_{1k}(t,\hat{u}'(t)) |u_k'|^2 u_k'^{*})+ a^{(1)}_2 2 Re(2 \hR^{(1)}_{2k}(t,\hat{u}'(t)) |u_k'|^2 u_k'^{*})$$ 
The equations for the rates of change of the $\hat{E}_i$ can be used for the computation of the $2 \times 2$ matrices $A$ and $B$ through the relations (\ref{rngmatrix1}) of Section \ref{algo}.

\section{The mesh refinement algorithm}\label{know2}
We continue our presentation with the mesh refinement algorithm. The construction in the previous section requires the exact knowledge of an accurate reduced model. This means, the knowledge of {\it both} the functional form of the reduced model and the associated coefficient vector $a^{(1)}.$ In fact, it is possible to relax this constraint by requiring the knowledge only of the functional form of the reduced model, i.e. knowledge of the vector $\hR^{(1)}$ but {\it not} of $a^{(1)}.$ This can be considered as a time-dependent generalization of the Swendsen renormalization algorithm (e.g. see the nice presentation in Ch. 5 of \cite{binney}), even though here we do not have a statistical framework. The Swendsen algorithm is based on the observation that knowledge of {\it only} the functional form of the reduced model but not necessarily of the associated coefficient vector $a^{(1)}$ is enough for computing quantities of the reduced system. In particular, the matrix $B$ can be calculated by using the resolved modes' values as computed from the full system.

As we have mentioned before, the entries of $B$ describe the contributions of the different terms appearing on the RHS of the reduced system to the rate of change of $E_i$ (the same for the entries of matrix $A$ and the full model). The determinant of the matrix $B$ measures whether there is need for the {\it reduced} system to transfer energy to smaller scales. The time instant when $\det B$ becomes nonzero, $T_B,$ signals the onset of energy transfer from the modes in $F$ to the modes in $G.$ The determinant of the matrix $A$ measures whether there is need for the {\it full} system to transfer energy to smaller scales. The time instant when $\det A$ becomes nonzero, $T_A,$ signals the onset of underresolution of the full system. The time interval $[T_B,T_A)$ is our window of opportunity to refine the mesh, without losing accuracy and without wasting computational resources.  We will use the value of $\det B$ as a criterion to decide when it is time to refine the mesh. 

Note that if there exists a singularity, the interval $\Delta T=T_A-T_B$ will shrink to zero as we increase the resolution.  The converse is not necessarily true. If $\Delta T$ appears to converge to zero as we increase the resolution does not mean that there certainly exists a singularity. Since all the calculations are finite, there is only a maximum resolution that we can afford. It may well be that an even larger, and presently unattainable, resolution could reveal that there is no singularity.   

The mesh refinement algorithm is given by:

\vskip14pt
{\bf Algorithm}
\begin{enumerate}
\item
Choose a value for $TOL.$ For this value of $TOL$ run a mesh refinement calculation, starting, say, from $N_{start}$ modes to $N_{final}$ modes. For example, let $N_{start}=32$ and double at each refinement until, say $N_{final}=256$ modes. Record the values of the quantities $\hat{E}_i, \; i=1,\ldots,m$ when $N=N_{final}$ and $|detB|=TOL.$ Let's call this simulation $S1.$
\item
For the same value of $TOL$ run a calculation with $N_{start}=N_{final}$ modes (for the example $N_{start}=N_{final}=256$). Record the values of the quantities $\hat{E}_i, \; i=1,\ldots,m$ when $|detB|=TOL.$ Let's call this simulation $S2.$
\item
Compare to within how many digits of accuracy the quantities $\hat{E}_i, \; i=1,\ldots,m$ computed from $S1$ and $S2$ agree. If the agreement is to within a specified accuracy, say 5 digits, then choose 
this value of $TOL.$ If the agreement is in fewer digits, then decrease $TOL$ (more stringent criterion) and repeat until agreement is met.
\item
Use the above decided value of $TOL$ to perform a mesh refinement calculation with a larger magnification ratio, i.e. a larger value for the ratio $N_{final}/ N_{start}.$  
\end{enumerate}
The agreement in digits of accuracy between $S1$ and $S2$ depends on the form of the terms chosen for the reduced model. Even though we do not know the coefficients of the reduced model, knowledge of the correct functional form of the terms can affect significantly the accuracy of the results. This situation is well known in the numerical study of critical exponents in equilibrium phase transitions (see e.g. Ch. 5 in \cite{binney}).


\subsection{How to compute the coefficients of the reduced model}\label{follow2}
When we only know the functional form of the terms appearing in the reduced model but not their coefficients it is not possible to evolve a reduced system. We present a way of actually computing the coefficients of the reduced model as needed. If the quantities $\hat{E}_i, \; i=1,\ldots,m$ are e.g. $L_p$ norms of the Fourier modes, then we can multiply Equations \eqref{reduced} with appropriate quantities and combine with Equations \eqref{conditions} to get
\begin{align*}\label{reduced}
\frac{d\hat{E}_1(\hat{u})}{dt} &= \sum_{i=1}^{m} a^{(1)}_i \hat{U}^{(1)}_{i1} (t,\hat{u}(t)) \\
\frac{d\hat{E}_2(\hat{u})}{dt} &= \sum_{i=1}^{m} a^{(1)}_i \hat{U}^{(1)}_{i2} (t,\hat{u}(t)) \\
 \quad     \cdots       \quad       & =   \quad       \cdots \quad \\
\frac{d\hat{E}_m(\hat{u})}{dt} &= \sum_{i=1}^{m} a^{(1)}_i \hat{U}^{(1)}_{im} (t,\hat{u}(t)) 
\end{align*}
where $\hat{U}^{(1)}_{ij}, \; i,j=1,\ldots,m$ are the new RHS functions that appear. Note that the RHS of the equations above does not involve primed quantities. The reason is that here the reduced quantities are computed by using the values of the resolved modes from the full system. The above system of equations is a linear system of equations for the vector of coefficients $a^{(1)}.$ In fact, the matrix of the system is the transpose $B^T$ of the matrix $B.$ The linear system can be written as
\begin{equation}\label{alphasystem}
B^{T} a^{(1)}= {\bf e}
\end{equation}
where ${\bf e}=\bigl(\frac{d\hat{E}_1(\hat{u})}{dt}, \ldots, \frac{d\hat{E}_m(\hat{u})}{dt} \bigr).$ This system of equations can provide us with the time evolution of the vector $a^{(1)}.$ 

The determination of coefficients for the reduced model through the system \eqref{alphasystem} is a time-dependent version of the method of moments. We specify the coefficients of the reduced model so that the reduced model reproduces the rates of change of a finite number of moments of the solution. This construction can actually be used as an adaptive way of determining a reduced model if one has access to experimental values of the rates of change of a finite number of moments. Suppose that we are conducting a real world experiment where we can compute the values of a finite number of moments on a coarse grid only. Then we can use the system \eqref{alphasystem} at predetermined instants to update a model defined on the coarse grid. Results of this construction will be presented elsewhere.


\section{Numerical results for the inviscid Burgers equation}\label{numerics2}

\begin{figure}
\centering
\subfigure[]{\epsfig{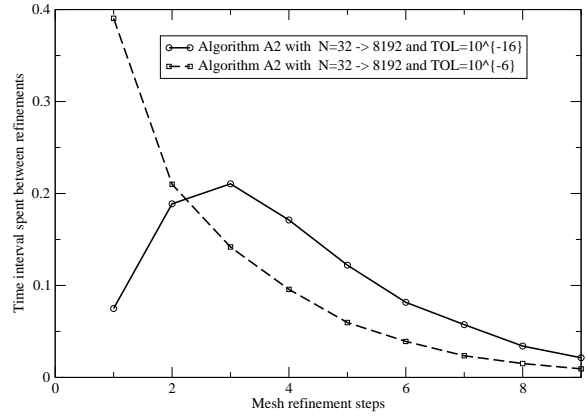}}
\qquad
\vskip20pt
\subfigure[]{\epsfig{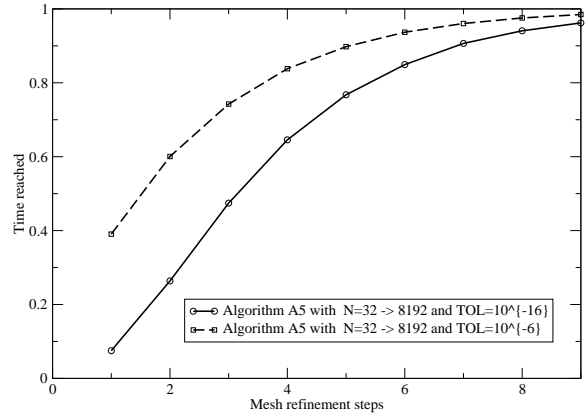}}
\caption{(a) Time spent between refinement steps for different tolerance values. (b) Time reached with the maximum allowed resolution. }
\label{plot_a5v0}
\end{figure}
                                                                                                                                                                                             
\begin{figure}
\centering
\epsfig{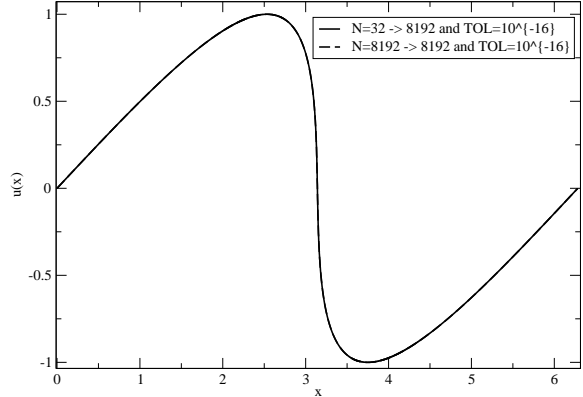}
\caption{Comparison of the velocity field produced at the time of termination of the mesh refinement algorithm for two different magnification ratios.
The first simulation has $N_{start}=32$ and $N_{final}=8192$ while the second has $N_{start}=N_{final}=8192.$ }
\label{plot_a5v0_velocity}
\end{figure}

We present numerical results of the mesh refinement algorithm for the inviscid Burgers equation with the initial condition $u(x,0)=\sin(x).$ This initial condition leads to a singularity forming at time $T_c=1.$ Figure \ref{plot_a5v0} contains results about the time spent between refinement steps and the time reached with the maximum allowed resolution. We start from a resolution $N_{start}=32$ and allow a maximum resolution of $N_{final}=8192.$ We present results for two values of the tolerance $TOL1=10^{-16}$ and $TOL2=10^{-6}.$ When the tolerance criterion is less strict the algorithm can reach later times before running out of resolution. 

In Figure \ref{plot_a5v0_velocity} we compare the velocity field produced by the algorithm with $N_{start}=32,$ $N_{final}=8192$ and $TOL1=10^{-16}$ with the velocity field produced by the algorithm with $N_{start}=N_{final}=8192$ and the same tolerance. It is obvious that the results are in very good agreement. However, the mesh refinement calculation was about 240 times faster. The final time reached by the algorithm is $T=0.962.$

\subsection{The direct approach to calculating the blow-up rate}
A mesh refinement algorithm can be used not only to approach a potential singularity but also estimate the rate at which the solution or some function of it blows-up. We restrict ourselves to the case of an algebraic (in time) singularity, meaning that some function of the solution diverges as $ \sim |T_c-T|^{-\gamma},$ where $\gamma > 0.$ Let us assume for a moment that $T_c$ is known. One obvious way of estimating $\gamma,$ is to run the mesh refinement algorithm and store the values of the blow-up quantity , say $\xi_n, \; n=1,\ldots,N,$  and the instant $T_n$ at which each refinement took place. Then, one can plot (in log-log) the values of the blow-up quantity at the different refinement instants $T_n$ as a function of the distance from the singularity $T_c-T_n$ and estimate the slope of the curve. That would provide us with the blow-up rate. Here we are interested in showing how the same estimate can be obtained using properties of a renormalization flow, i.e. a coarse-graining (dimensional reduction) procedure. Before we proceed, we have to address the issue of the value of $T_c$ which is, in general, unknown. Thus, the value of $T_c$ has to be calculated from the algorithm. It is simple to see that small errors in the estimation of $T_c$ can lead to huge errors in the estimation of the blow-up rate. One way of estimating $T_c$ is the following : for different choices of $T_c,$ plot, in log-log coordinates, the values of the blow-up quantity at the refinement instants $T_n$ as a function of the distance from the singularity $T_c-T_n$  and pick the value of $T_c$ for which this plot is a straight line. This can be decided by monitoring the value of the correlation coefficient for a linear regression. 

\begin{figure}
\centering
\epsfig{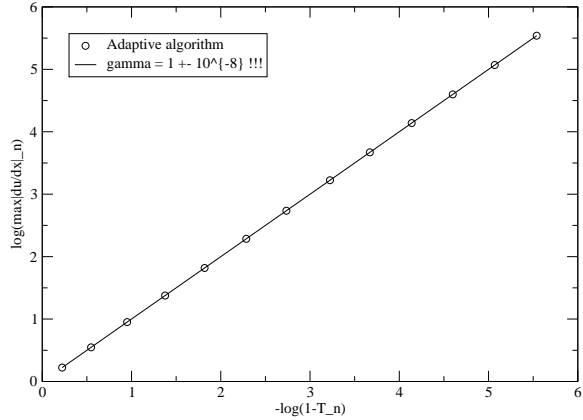}
\caption{Log-log plot of the maximum absolute value of the velocity gradient ${\max |\p{u}{x}|}_n$ and  $(1-T_n)^{-1}$ for the different refinement steps (indexed by $n$).}
\label{plot_maxvor_time}
\end{figure}

We present results of the above construction for the inviscid Burgers equation with the initial condition $u(x,0)=\sin(x).$ This initial condition leads to a singularity forming at time $T_c=1.$ Figure \ref{plot_maxvor_time} shows the log-log plot of the maximum absolute value of the velocity gradient $\log({\max |\p{u}{x}|}_n)$ and of the inverse distance from the singularity time $(1-T_n)^{-1}$ as recorded at the different refinement steps $T_n$. The slope of the curve is $\gamma = 1 \pm 10^{-8}.$ Note that the minute error in the estimate shows that the refinement algorithm keeps the calculation well-resolved even very close to the singularity. The calculations were performed using the mesh refinement algorithm of Section \ref{know2} with the refinement tolerance criterion $TOL=\det B$ set to $10^{-10}.$ We should note that for this value of $TOL,$ the value of $\det A$ for the full system is still much smaller than the double precision roundoff threshold of $10^{-16}.$ For this calculation we set $N_{start}=32$ and $N_{final}=131072$ and the algorithm terminated at time $T=0.996.$ The mesh refinement is about 3000 times faster than a calculation with $N_{start}=N_{final}=131072.$

\section{Numerical results for the supercritical focusing Schr\"odinger equation}\label{schrodinger}

\begin{figure}
\centering
\includegraphics[width=3.5in]{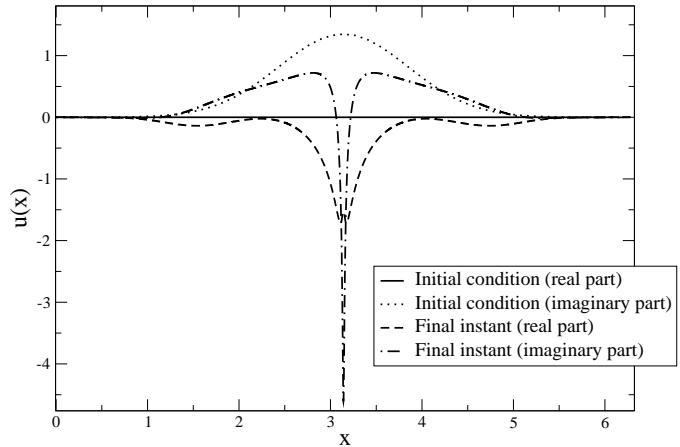}
\caption{Supercritical ($\sigma=3$) Schr\"odinger equation with $\max |u_0(x)|=1.35.$}
\label{supercritical_1}
\end{figure}

We continue with numerical results about the supercritical focusing Schr\"odinger equation. The focusing Schr\"odinger equation is given by 
\begin{equation}
i \p{u}{t} + \Delta u + |u|^{2\sigma}u=0, \, \, \text{where} \, \, \sigma >0
\label{schrod} 
\end{equation}
The equation needs to be supplemented by an initial condition $u(x,0)=u_0(x)$ and boundary conditions. It has been conjectured by Zakharov \cite{zakharov} that in $d$ dimensions, when $\sigma > 2/d$ and for sufficiently large initial condition, the solution of \eqref{schrod} will blow-up at a finite time $T$ and the behavior of the solution close to the blow-up time is given by 
$$u(x,t) = ((2\kappa(T-t))^{-\frac{1}{2}(\frac{1}{\sigma}+i\frac{\omega}{\kappa})}Q((2\kappa(T-t))^{-1/2}|x|),$$
where $Q(\xi)$ is a complex-valued function with appropriate decay properties and $\kappa$ and $\omega$ are parameters to be determined. For the maximum of the solution we have
$$ \max |u(x,t)| \sim (T-t)^{-\frac{1}{2\sigma}}  \, \, \text{as} \, \, t \rightarrow T.$$
Although the mathematical theory is not yet complete, overwhelming evidence from numerical and formal analytical calculations suggests that the conjecture is true. Here, we restrict attention to the 1D case and to periodic boundary conditions in the domain $[0,2\pi].$ In the 1D case, according to the conjecture the solution exhibits a algebraic finite time blow-up when $\sigma > 2.$ Here we present results for the case $\sigma=3.$ In the numerical experiments we used the initial condition 
$$ u_0(x,0)=i A \exp(-(x-\pi)^2),$$
for different values of $A.$ For this initial condition we have $\max |u_0(x)| = A$ at $x=\pi.$ Figure \ref{supercritical_1} shows the initial condition for $A=1.35$ and the solution as computed by the mesh refinement algorithm with $N_{start}=48$ and $N_{final}=10368.$ The tolerance criterion $TOL=\det B$ was set to $10^{-16}.$ The algorithm was implemented with the $t$-model for the reduced system as in the case of inviscid Burgers.

Table \ref{supercritical_2} contains the estimated blow-up exponents for the maximum of the solution for different values of $A.$ For $A=1.242$ the mesh refinement algorithm does not run out of resolution which signals the absence of a singularity. For all the other cases and for $N_{start}=48$ and $N_{final}=10368,$ the mesh refinement algorithm was about 200 times faster than a calculation performed with $N_{start}=N_{final}=10368.$ Unlike the case of the inviscid Burgers equation, here we cannot estimate beforehand the exact time $T$ of the blow-up. We do that in the way proposed in the previous section. In particular, for different choices of $T,$ we calculated the correlation coefficient of the linear fit (in log-log coordinates), of the values of the blow-up quantity as a function of the distance from the singularity $T-t$  and picked the value of $T$ for which the correlation coefficient is largest . For all the cases shown in Table \ref{supercritical_2} the correlation coefficient is about 0.999999999. The algorithm is able to approach the estimated singularity instant $T$ to within $5\times 10^{-5}$ units of time. The conjectured blow-up exponent for the maximum of the solution when $\sigma=3$ is $1/2\sigma=1/6 \sim 0.1667.$ The relative deviation of the estimated values of the exponent relative to the conjectured value of the exponent is within $1\%$ for all the values of $A$ examined except for $A=1.8$ and $A=2.$ 

We would like to make a comment about the discrepancy for $A=1.8$ and $A=2.$. It is to be expected that if one keeps the same maximum resolution while increasing the magnitude of the initial condition, i.e. the value of $A$, after some value of $A$ the algorithm runs out of resolution before it can come close enough to the singularity for the asymptotic behavior to settle in. To elucidate this point we also ran the mesh refinement algorithm with $N_{final}=34992$ for $A=1.8$ and $A=2.$ The estimated values of the blow-up exponent are included in Table \ref{supercritical_2} in parentheses. As we see, if one uses large enough resolution, the relative deviation of the estimated values of the exponent relative to the conjectured value of the exponent decreases again to within $1\%.$ Note that for $N_{start}=48$ and $N_{final}=34992$ the mesh refinement algorithm is about 400 times faster than a calculation with $N_{start}=N_{final}=34992.$  As expected, the acceleration factor increases when $N_{final}/N_{start}$ increases.

\begin{table}
\centering
\begin{tabular}{| c | c | c | c |}
\hline
$\max |u_0(x)|$ & Est. exp. $\alpha$ & $| \text{Rel. dev.} |$ ($\%$) & Sing. form \\  \hline
1.242   & ---        &  ---        & ---\\ \hline
1.243   & 0.1652  & 0.90  &   $(T-t)^{-\alpha}$ \\ \hline
1.250   & 0.1648  & 1.14  &\\ \hline
1.255   & 0.1654  & 0.78  &\\ \hline
1.260   & 0.1649  & 1.08  &\\ \hline
1.300   & 0.1655  & 0.72  &\\ \hline
1.350   & 0.1662  & 0.30  &\\ \hline
1.500   & 0.1678  & 0.66  &\\ \hline
1.600   & 0.1691  & 1.44  & \\ \hline
1.800   & 0.1727 (0.1684)  & 3.60 (1.01)  &\\ \hline
2.000   & 0.1766 (0.1696)  & 5.93 (1.74)  &\\ \hline
\end{tabular}
\caption{Estimated blow-up exponents for the supercritical ($\sigma=3$) Schr\"odinger equation.  The relative deviation is from the conjectured value of $1/2\sigma=1/6 \sim 0.1667.$} 
\label{supercritical_2}
\end{table}


\section{Calculation of the blow-up rate as a critical exponent}\label{rate}

As we have said, we are also interested in showing how the blow-up rate estimate can be obtained using properties of a renormalization flow, i.e. a coarse-graining process. There are two ways to do that:  i) Wilson's or ``phase transition" approach, where one treats the singularity as a fixed point of a renormalization transformation and computes the blow-up rate by analysis in the vicinity of the fixed point, and ii) the Widom-Kadanoff or ``scaling approach", where one assumes the existence of certain scaling laws in the vicinity of the singularity and then combines them to obtain the blow-up rate. 

\subsection{The ``phase transition" approach}\label{phase}
The key idea is  that a series of successive refinement steps (going to smaller and smaller scales) can be seen (approximately) as a coarse-graining process in reverse. Thus, one can run the mesh refinement algorithm, compute and store the coefficients of the reduced model at each refinement step and then use them to reconstruct the renormalization flow from smaller to larger scales. In this case, the smallest scale that the refinement algorithm reached is the starting scale of the renormalization flow. For the case of a time-dependent PDE the mesh refinement algorithm allows us to get closer and closer to the singularity instant $T_c$. Thus, the renormalization procedure will take us further and further away from $T_c.$

There are two ways to show how the renormalization flow can be used to compute the blow-up rate. The first, the ``phase transition" approach, assumes that the phase transition is a fixed point of the renormalization flow and proceeds with an analysis near the fixed point (e.g. pp. 124-127 in \cite{binney}). However, as we mentioned at the discussion after \eqref{rngmatrix2}, we do not use a linear stability analysis because the eigenvalues of $M$ vary most rapidly near the fixed point. Instead, we deal with the full (nonlinear) renormalization flow.

The second way, the ``scaling" approach, is just a manipulation of different scaling laws assuming to hold asymptotically near the singularity. Of course, both lead to the same expression for the blow-up rate. We choose to present both since it elucidates further the connection between the techniques presented in this paper and those used in the theory of equilibrium phase transitions.

\begin{figure}
\centering
\epsfig{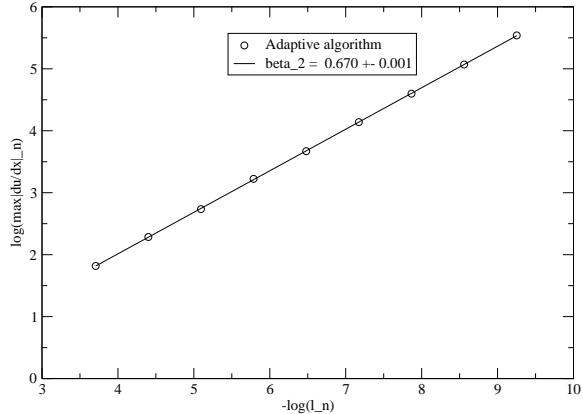}
\caption{Log-log plot of the maximum absolute value of the velocity gradient ${\max |\p{u}{x}|}_n$ and the inverse length scale of the reduced system $l_n^{-1}$ for the different renormalization steps (indexed by $n$).}
\label{plot_maxvor_scale}
\end{figure}

\begin{figure}
\centering
\epsfig{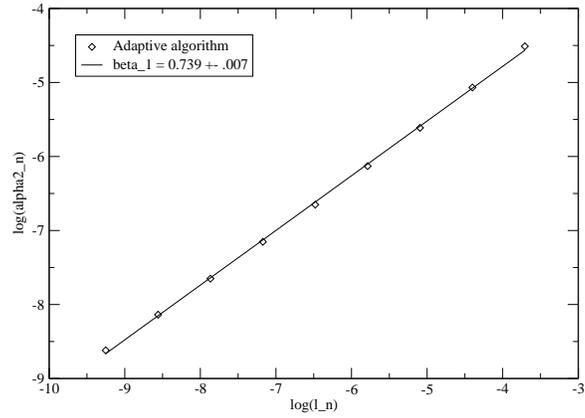}
\caption{Log-log plot of the coefficient $a_2^{(n)}$ of the $t$-model and the length scale of the reduced system $l_n$ for the different renormalization steps (indexed by $n$).}
\label{plot_alpha2_scale}
\end{figure}

\begin{figure}
\centering
\epsfig{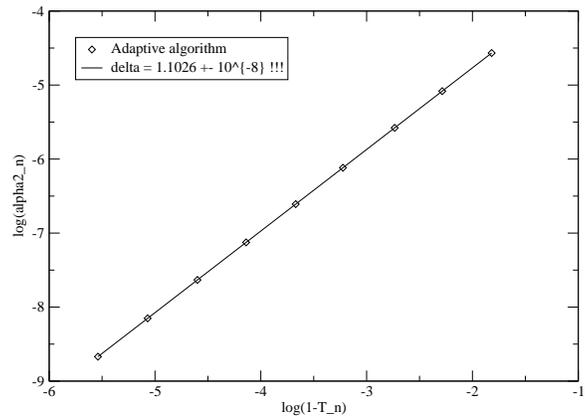}
\caption{Log-log plot of the coefficient $a_2^{(n)}$ of the $t$-model  and  $1-T_n$ for the different renormalization steps (indexed by $n$).}
\label{plot_alpha2_time}
\end{figure}

We start our presentation of the blow-up rate calculation with the ``phase transition" approach (see e.g. \cite{binney}). Let us suppose that near the singularity instant $T_c$ a quantity $\xi$ behaves as 
$ |T_c-T|^{-\gamma}.$ For the case of Burgers this would be the maximum of the velocity gradient, i.e. 
$\max |\p{u}{x}|.$ We want to find the value of $\gamma.$ As we have said we assume that we have computed and stored a sequence of coefficients for the reduced model, the associated length scale, the value of the blow-up quantity and the time of occurrence of the refinement step. Then, by simply reversing the sequence indexing, we have the necessary quantities for the description of a renormalization flow which starts close to $T_c$ and moves further away with every coarse-graining step. Since every renormalization step brings us further away from the critical point $T_c$, the values of the blow-up quantity become smaller with every renormalization step. Thus, if we coarse-grain the length scale at which we probe the problem by a factor of $b$ at each step (where $b > 1$), then 
$\xi_{n+1}=\frac{\xi_n}{b^{\beta_2}},$ with $\beta_2 > 0.$ This implies $\xi_n \sim l_n^{-\beta_2}$ and thus $\beta_2$ can be computed from the refinement algorithm data collected. The coefficient of the reduced model which monitors the deviation of the full and reduced model will increase with each renormalization step, i.e. $\alpha_{n+1}=\alpha_n b^{\beta_1}, \; \text{with } \beta_1 > 0.$ This implies 
$\alpha_n \sim l_n^{\beta_1}$ and $\beta_1$ can also be computed from the collected data. Moreover, repeated application of the recursive relation for the coefficient $\alpha_{n}$ gives 
$\alpha_n = \alpha_0 (b^{\beta_1})^n.$ This relation is the analog of the recursive relation derived in the theory of phase transitions by linearization of the renormalization flow around the critical (fixed) point. Here we did {\it not} resort to a linearization procedure. To proceed, we need to estimate the behavior of $\alpha_0,$ the starting point of the renormalization flow. In the theory of phase transitions, the behavior of the coefficient $\alpha_0$ is assumed to be linear in $|T_c-T|.$ However, there is no {\it a priori} reason for such a behavior. We assume that 
$\alpha_0=C_2  |T_c-T|^{\delta},$ where $\delta$ can also be computed from the collected data. 

Let us summarize what we have obtained so far. As we renormalize, the blow-up quantity decreases and the reduced model coefficient that monitors the deviation of the full and reduced model increases. Following the phase transition approach we thus assume that if we take enough renormalization steps then we have 
$$ \frac{\xi}{C_1(b^{\beta_2})^n}=u   \; \; \text{and} \; \; C_2 |T_c-T|^{\delta}(b^{\beta_1})^n= v$$
where $u,v$ are quantities of the same order and $C_1,C_2$ are constants that depend on the initial conditions. We can eliminate $n$ in the above two relations and get 
$$ \xi \sim  |T_c-T|^{-\gamma}, \; \text{with} \; \gamma=\frac{\delta \beta_2}{\beta_1}.$$
Thus, we have expressed the blow-up rate exponent $\gamma$ as a function of scaling exponents that are associated with properties of the renormalization flow. 

Before we conclude with this approach, we need to make one more comment. We have said before that the ``phase transition" approach treats the singularity as a fixed point of the renormalization flow. To do that one has to construct a differential equation for the evolution of the coefficient $\alpha$ with respect to $l.$ Note that by the way we have defined it, $\alpha$ is dimensionless. The equation for its evolution with changes in $l$ is given by $l \p{\alpha}{l}=\beta(\alpha)$ \cite{binney}. The RHS of the equation is called the beta function and its zeros determine the fixed points of the renormalization flow. Since $\alpha =C l^{\beta_1},$ for some constant $C,$ we have $  l \p{\alpha}{l} =C\beta_1 l^{\beta_1}.$ So, the beta function is $\beta(\alpha)=C \beta_1 l^{\beta_1}=\beta_1 \alpha.$ So, the only fixed point of the beta function is $\alpha=0.$ If $\beta_1 >0$ then $\alpha=0$ corresponds to $l=0,$ i.e. the zero scale. But this is exactly the active scale reached at the instant that the singularity occurs. So, the singularity is indeed a fixed point of the renormalization flow as long as $\beta_1 >0$. Moreover, if $\beta_1 > 0$, this fixed point is unstable, so that if we start close to it, the renormalization flow will take us further away. This is indeed the case for the Burgers equation as we show numerically in the next section.  

This concludes the ``phase transition" approach.

\subsection{The ``scaling" approach}\label{scaling}
We conclude with the ``scaling" approach which is based on direct combination of the different scaling laws associated with the renormalization flow. Indeed, let $  \xi \sim  |T_c-T|^{-\gamma'},$ where $\gamma'$ is the blow-up rate exponent to be estimated. If we assume that near $T_c$ we have $\xi \sim l^{-\beta_2},$ $\alpha \sim l^{\beta_1}$ and $\alpha \sim |T_c-T|^{\delta},$ we can use the renormalization flow to estimate $\beta_1,$ $\beta_2$ and $\delta.$ Then a straightforward combination of the three scaling laws leads to $\gamma'=\frac{\delta \beta_2}{\beta_1}.$ So, $\gamma'=\gamma$ and as expected this approach leads to the same expression for the blow-up rate exponent as the "phase transition" approach.

Figures \ref{plot_maxvor_scale}-\ref{plot_alpha2_time} show how one can use the above construction to estimate the blow-up rate $\gamma$ from renormalization flow quantities. Recall that the coefficient of the reduced model that monitors the deviation of the reduced and full systems is $a_2^{(n)}.$ Also, that the index $n$ appearing in the figures is used now to count the renormalization steps which are the opposite of the refinement steps.  The length scale $l_n$ at which we probe the system for the different renormalization steps is the length scale of the reduced model. This means that if we have a full system calculation with $N_n$ modes, then $l_n=2\frac{2\pi}{N_n},$ since the reduced model has half the resolution of the full system. 

From the data we estimate the exponents $\beta_2=0.670 \pm 0.001,$ $\beta_1=0.739 \pm 0.007,$ and $\delta = 1.1026 \pm 10^{-9}.$ From these estimates we get $\gamma' = 1 \pm 0.01.$ Thus, when we compute the blow-up rate using solely renormalization flow quantities, the estimation error is larger than when computing this rate directly. This is to be expected since we had to combine three empirically determined scaling laws, each one of which comes with its own error and also relies entirely on the adequacy of the reduced model. Nevertheless, the obtained accuracy is acceptable and moreover, it highlights the accuracy of the $t$-model for this equation.  

Finally, since $\beta_1=0.739 > 0,$ we conclude that the singularity is an unstable fixed point of the renormalization flow (see discussion at the end of Section \ref{phase}).
 

\section{Conclusions and future work}\label{conclusions}
We have presented a mesh refinement algorithm, inspired by renormalization constructions in critical phenomena, which allows the efficient location and approach of a possible singularity. The algorithm assumes knowledge of an accurate reduced model. In particular, it assumes knowledge of the functional form of the reduced model but not of the actual coefficients. We provide a way of computing the necessary coefficients on the fly as needed. On a theoretical level, the algorithm can be used to study the behavior of (near-) singular solutions. On the practical side, it can be used as a mesh refinement tool.

We have only examined the simple case of periodic boundary conditions and the mesh refinement performed was uniform. We plan to extend the constructions presented here to a real space formulation which will allow the treatment of non-periodic boundary conditions and more complicated geometries. In that case, one can divide the domain into sub-domains and apply the mesh refinement algorithm individually in the different sub-domains. In addition, the algorithm can be modified to perform mesh-coarsening after the computationally intensive time interval of the simulation has passed.

The original motivation behind the development of the algorithm was the open problem of the formation of singularities in finite time for the incompressible Euler and Navier-Stokes equations of fluid mechanics. In addition to helping with the issue of singularity formation, we hope that the algorithm can be of use in the simulation of real world flows by allowing a better assessment of the onset of underresolution.

\section*{Acknowledgements} I am grateful to Profs. G.I. Barenblatt, 
A.J. Chorin and O.H. Hald for their ongoing guidance and support. I would like to thank Prof. V. Sverak for helpful discussions and Profs. S. Weinberg and K. Wilson for inspiration.

\end{document}